# On The Existence of Periodic Solutions for a Certain System of Third Order Nonlinear Differential Equations

MUZAFFER ATES
Van Technical School , University of Yuzuncu Yıl, 65080, Van, TURKEY
ates.muzaffer@yahoo.com

## Abstract

In this paper, we study the existence and uniqueness of periodic solutions of the differential equation of the form

$$\dddot{X} + F(X, \dot{X}, \ddot{X})\ddot{X} + G(X, \dot{X})\dot{X} + H(X) = P(t, X, \dot{X}, \ddot{X}).$$

Here, we obtain some sufficient conditions which guarantee the existence of periodic solutions. This equation is a quite general third- order nonlinear vector differential equation, and one example is given for illustration of the subject.

## 1. Introduction

There have been done many studies concerning the problem of qualitative behaviors of solutions of certain third order nonlinear scalar and vector differential equations, see $[1-11]$. However, there are only a few papers on the existence and uniqueness of periodic solutions of third order nonlinear vector differential equations without any example. Some of them can be summarized here as follows:

In 1995, Feng [3] considered the differential equation of the form

$$\dddot{X} + A(t)\ddot{X} + B(t)\dot{X} + H(X) = P(t, X, \dot{X}, \ddot{X}).$$

He proved the existence and uniqueness of periodic solution. Later, Tiryaki [6] investigated the boundedness and periodicity results of the solutions of vector differential equation

$$\dddot{X} + A\ddot{X} + G(\dot{X}) + H(X) = P(t, X, \dot{X}, \ddot{X}).$$

Similarly, Tunç [7] proved some results on the boundedness and periodicity of the solutions of the vector differential equation

$$\dddot{X} + F(X, \dot{X})\ddot{X} + B\dot{X} + H(X) = P(t, X, \dot{X}, \ddot{X}).$$

Recently, Tunç and Ateş [9] studied the existence and uniqueness of periodic solutions of third order nonlinear differential equations

$$\dddot{X} + A(t)\ddot{X} + G(\dot{X}) + H(X) = P(t, X, \dot{X}, \ddot{X}),$$

and

$$\dddot{X} + F(X, \dot{X})\ddot{X} + B(t)\dot{X} + H(X) = P(t, X, \dot{X}, \ddot{X}).$$

In this paper, we consider the nonlinear vector differential equation

$$\dddot{X} + F(X, \dot{X}, \ddot{X})\ddot{X} + G(X, \dot{X})\dot{X} + H(X) = P(t, X, \dot{X}, \ddot{X}) \tag{1.1}$$

where $X \in R^n$ and $t \in [0, \infty)$; $F$ and $G$ are $n \times n$ - symmetric continuous matrix functions; $H : R^n \to R^n$ and $P : R^+ \times R^n \times R^n \times R^n \to R^n$, and $P$ is a periodic fuction, that is,

$$P(t + \omega, X, \dot{X}, \ddot{X}) = P(t, X, \dot{X}, \ddot{X}), \quad \omega > 0 \text{ is period.}$$

Given any $X, Y$ in $R^n$, the symbol $\langle X, Y \rangle$ is used to denote the usual scalar product in $R^n$, that is, $\langle X, Y \rangle = \sum_{i=1}^{n} x_i y_i$, thus $\langle X, X \rangle = \|X\|^2$.

Throughout this paper we assume that the following:

There exist $n \times n$ real constant symmetric matrices $A$, $B$ and an $n \times n$ operator $A(X, Y)$, such that

$$H(X) = H(Y) + A(X, Y)(X - Y) \tag{1.2}$$

for which the eigenvalues $\lambda_i(A(X, Y))$ are continuous and satisfy

$$0 < \delta_h \leq \lambda_i(A(X, Y)) \leq \Delta_h \tag{1.3}$$

for fixed constants $\delta_h$ and $\Delta_h$.

We shall assume that $\Delta_h \leq k \delta_a \delta_b$, $(k < 1)$

where $$k = \min \frac{1}{8} \left\{ \frac{1}{2}, \frac{\delta_b}{\delta_a \Delta_a} \right\}. \tag{1.4}$$

The eigenvalues of the related matrices are such that

$$0 < \delta_a = \min\{\lambda_i(A), \lambda_i(F(X, Y, Z))\}, \quad \Delta_a = \max\{\lambda_i(A), \lambda_i(F(X, Y, Z))\}$$

$$0 < \delta_b = \min\{\lambda_i(B), \lambda_i(G(X, Y))\}, \quad \Delta_b = \max\{\lambda_i(B), \lambda_i(G(X, Y))\}$$

and

$$0 < \lambda_i(F(X, Y, Z) - A) \leq \frac{\sqrt{\varepsilon}}{2}, \quad 0 < \lambda_i(G(X, Y) - B) \leq \frac{\sqrt{\varepsilon}}{2}$$
$$(i = 1, 2, \ldots, n),$$

where

$$\sqrt{\varepsilon} \leq \min \left\{ \frac{\delta_b \delta_h}{4\Delta_b + 4}, \frac{\delta_a \delta_b}{6\Delta_a + 7}, \frac{\delta_a}{2}, 1 \right\}. \tag{1.5}$$

**Remark.** Motivation of this study has been based on that of Feng [3], Tiryaki [6], Tunç [7] Tunç and Ateş [9]. Equation (1.1) is a quite general third- order nonlinear vector differential equation. In particular, many third-order differential equations which have been discussed in [1-11] are special cases of Eq. (1.1).

## 2. Main Result

**Theorem:** Suppose that
(i) there exists an $n \times n$ - real continuous operator $A(X,Y)$ for any vectors $X,Y$ in $R^n$, such that
$$H(X) = H(Y) + A(X,Y)(X - Y) \tag{2.1}$$
whose eigenvalues $\lambda_i(A(X,Y))$ $(i = 1,2,...,n)$ satisfy
$$0 < \delta_h \leq \lambda_i(A(X,Y)) \leq \Delta_h$$
for fixed constants $\delta_h$ and $\Delta_h$, and
$$\Delta_h \leq k\delta_a \delta_b$$
where the positive constant $k$ to be determined later in the proof;

(ii) the symmetric matrices $F$ and $G$ have positive eigenvalues and commute with themselves as well as with the operator $A(X,Y)$ for any vector $X,Y,Z$ in $R^n$, and $X,Y$ in $R^n$, respectively;

(iii) there exist finite constants $\delta_0 \geq 0, \delta_1 \geq 0$ such that the vector $P$ satisfies

$$\|P(t,X,Y,Z)\| \leq \delta_0 + \delta_1(\|X\| + \|Y\| + \|Z\|) \tag{2.2}$$

uniformly in $t \geq 0$ for all arbitrary $X, Y, Z$ in $R^n$;

(iv) let $0 < \varepsilon \leq 1$

where
$$\sqrt{\varepsilon} \leq \min\left\{\frac{\delta_b \delta_h}{4\Delta_b + 4}, \frac{\delta_a \delta_b}{6\Delta_a + 7}, \frac{\delta_a}{2}, 1\right\}. \tag{2.3}$$

Then, if $H(0) = 0$ and $\delta_1$ is sufficiently small, then Eq. (1.1) has at least a periodic solution.

If $P(t,X,Y,Z) = P(t)$, Eq. (1.1) has a unique periodic solution. Then, the condition (2.2) can be improved to

$$\|P(t,X,Y,Z)\| \leq \theta_1(t) + \theta_2(t)(\|X\|^2 + \|Y\|^2 + \|Z\|^2)^{\frac{1}{2}} \tag{2.4}$$

where $\theta_1(t)$ and $\theta_2(t)$ are continuous functions of $t$ satisfying
$$0 \leq \theta_1(t) < \alpha_0,$$
$$\text{for all } t \text{ in } R. \tag{2.5}$$
$$0 \leq \theta_2(t) < \alpha_1$$

In the subsequent discussion we require the following lemmas.
**Lemma 1:** Let $D$ be a real symmetric $n \times n$ matrix, then for any $X$ in $R^n$ we have

$$\delta_d \|X\|^2 \leq \langle DX, X \rangle \leq \Delta_d \|X\|^2$$

where $\delta_d$, $\Delta_d$ are the least and the greatest eigenvalues of $D$, respectively.
**Proof:** See [11].

**Lemma 2**: Let $Q$, $D$ be any two real $n \times n$ commuting symmetric matrices.

Then
(i) the eigenvalues $\lambda_i(QD)$ ($i = 1,2,...,n$) of the product matrix $QD$ are all real and satisfy

$$\max_{1 \leq j,k \leq n} \lambda_j(Q)\lambda_k(D) \geq \lambda_i(QD) \geq \min_{1 \leq j,k \leq n} \lambda_j(Q)\lambda_k(D);$$

(ii) the eigenvalues $\lambda_i(Q + D)$ ($i = 1,2,...,n$) of the sum of the matrices $Q$ and $D$ are all real and satisfy

$$\left\{\max_{1 \leq j \leq n} \lambda_j(Q) + \max_{1 \leq k \leq n} \lambda_k(D)\right\} \geq \lambda_i(Q+D) \geq \left\{\min_{1 \leq j \leq n} \lambda_j(Q) + \min_{1 \leq k \leq n} \lambda_k(D)\right\}.$$

**Proof:** See [11].

## 3. Proof of the Theorem

**Proof.** Our main toll in the proof is the vector Lyapunov function

$V = V(t, X, Y, Z)$ defined by

$$2V = \frac{1}{4}\langle BX, BX \rangle + \frac{3}{2}\langle BY, Y \rangle + \langle Z, Z \rangle + \langle Z + AY + \frac{1}{2}BX, Z + AY + \frac{1}{2}BX \rangle \tag{3.1}$$

where $A$ and $B$ are real $n \times n$ constant symmetric matrices.
Then, there exist positive constants $\delta_2$ and $\delta_3$ such that

$$\delta_2\left(\|X\|^2 + \|Y\|^2 + \|Z\|^2\right) \leq 2V \leq \delta_3\left(\|X\|^2 + \|Y\|^2 + \|Z\|^2\right). \tag{3.2}$$

Let us, for convenience, replace Eq. (1.1) by the equivalent form

$$\begin{cases} \dot{X} = Y, \dot{Y} = Z \\ \dot{Z} = -F(X,Y,Z)Z - G(X,Y)Y - H(X) + P(t,X,Y,Z) \end{cases} \tag{3.3}$$

Let $(X, Y, Z)$ be any solution of (3.3), then the total derivative of $V$ with respect to $t$ along this solution path is

$$\dot{V} = \frac{d}{dt}V[X(t), Y(t), Z(t)] = -V_1 - V_2 - V_3 + V_4 \tag{3.4}$$

where

$$V_1 = \frac{1}{8}\langle BX, H(X) \rangle + \langle H(X), AY \rangle + \frac{1}{4}\langle AY, G(X,Y)Y \rangle$$

$$V_2 = \frac{1}{8}\langle BX, H(X)\rangle + \frac{1}{2}\langle F(X,Y,Z)Z, Z\rangle + 2\langle H(X), Z\rangle$$

$$V_3 = \frac{1}{4}\langle BX, H(X)\rangle + \frac{1}{4}\langle AY, G(X,Y)Y\rangle + \frac{1}{2}\langle F(X,Y,Z)Z, Z\rangle$$

$$+ \frac{1}{2}\langle BX, (F(X,Y,Z) - A)Z\rangle + \frac{1}{2}\langle BX, (G(X,Y) - B)Y\rangle$$

$$+ \langle AY, (F(X,Y,Z) - A)Z\rangle + 2\langle (G(X,Y) - B)Y, Z\rangle$$

$$+ \langle (F(X,Y,Z) - A)Z, Z\rangle + \frac{1}{2}\langle (G(X,Y) - B)Y, AY\rangle$$

$$V_4 = \langle \frac{1}{2}BX + AY + 2Z, P(t, X, Y, Z)\rangle .$$

From (1.2) we have
$$H(X) = H(0) + A(X, 0)X .$$
Thus, if $H(0) = 0$ and condition (1.3) is satisfied, we obtain the following inequalities

$$\langle BX, H(X)\rangle = \langle BX, A(X,0)X\rangle \geq \delta_b \delta_h \|X\|^2 ;$$

$$\langle AY, G(X,Y)\rangle \geq \delta_a \delta_b \|Y\|^2 ;$$

$$\langle F(X,Y,Z)Z, Z\rangle \geq \delta_a \|Z\|^2 .$$

Next, we give estimates for the other terms of $\dot{V}$.

For some constants $k_j > 0$, $(j = 1, 2, ..., 6)$, conveniently chosen later, we obtain
$$\langle H(X), AY\rangle = \frac{1}{2}\|k_1^{-1}(H(X) + k_1 AY)\|^2 - \frac{1}{2}k_1^{-2}\langle H(X), H(X)\rangle - \frac{1}{2}k_1^2\langle AY, AY\rangle$$

$$\geq -\frac{1}{2}k_1^{-2}\delta_h \Delta_h \|X\|^2 - \frac{1}{2}k_1^2 \delta_a \Delta_a \|Y\|^2 ;$$

in a similar way we have the following

$$2\langle H(X), Z\rangle \geq -k_2^{-2}\delta_h \Delta_h \|X\|^2 - k_2^2 \|Z\|^2 ;$$

$$\frac{1}{2}\langle BX, (F(X,Y,Z) - A)Z\rangle = \frac{1}{4}\|k_3^{-1}\sqrt{B}\sqrt{F - A}X + k_3\sqrt{B}\sqrt{F - A}Z\|^2$$

$$-\frac{1}{4}k_3^{-2}\langle BX,(F-A)X\rangle - \frac{1}{4}k_3^2\langle BZ,(F-A)Z\rangle$$

$$\geq -\frac{1}{8}k_3^{-2}\Delta_b\sqrt{\varepsilon}\|X\|^2 - \frac{1}{8}k_3^2\Delta_b\sqrt{\varepsilon}\|Z\|^2$$

$$\geq -\Delta_b\sqrt{\varepsilon}\|X\|^2 - \frac{1}{3}\sqrt{\varepsilon}\|Z\|^2 \text{ for } k_3^2 = \min\left\{\frac{1}{8},\frac{8}{3\Delta_b}\right\};$$

$$\frac{1}{2}\langle BX,(G(X,Y)-B)Y\rangle \geq -\frac{1}{8}k_4^{-2}\Delta_b\sqrt{\varepsilon}\|X\|^2 - \frac{1}{8}k_4^2\Delta_b\sqrt{\varepsilon}\|Y\|^2$$

$$\geq -\sqrt{\varepsilon}\|X\|^2 - \frac{7}{4}\sqrt{\varepsilon}\|Y\|^2 \text{ for } k_4^2 = \min\left\{\frac{\Delta_b}{8},\frac{14}{\Delta_b}\right\};$$

$$\langle AY,(F(X,Y,Z)-A)Z\rangle \geq -\frac{1}{4}k_5^{-2}\Delta_a\sqrt{\varepsilon}\|Y\|^2 - \frac{1}{4}k_5^2\Delta_a\sqrt{\varepsilon}\|Z\|^2$$

$$\geq -\frac{3}{4}\Delta_a\sqrt{\varepsilon}\|Y\|^2 - \frac{1}{3}\sqrt{\varepsilon}\|Z\|^2 \text{ for } k_5^2 = \min\left\{\frac{1}{3},\frac{4}{3\Delta_a}\right\};$$

$$2\langle Z,(G(X,Y)-B)Y\rangle \geq -k_6^{-2}\frac{\sqrt{\varepsilon}}{2}\|Y\|^2 - k_6^2\frac{\sqrt{\varepsilon}}{2}\|Z\|^2$$

$$\geq -\frac{3}{4}\Delta_a\sqrt{\varepsilon}\|Y\|^2 - \frac{1}{3}\sqrt{\varepsilon}\|Z\|^2 \text{ for } k_6^2 = \min\left\{\frac{2}{3\Delta_a},\frac{2}{3}\right\};$$

and we are left with

$$\langle (F(X,Y,Z)-A)Z,Z\rangle + \frac{1}{2}\langle (G(X,Y)-B)Y,AY\rangle \geq 0$$

because

$$\lambda_i[F(X,Y,Z)-A]\|Z\|^2 \geq 0, \qquad \lambda_i(A)\lambda_i[G(X,Y)-B]\|Y\|^2 \geq 0.$$

Then, rearranging the terms of $V_1, V_2$ and $V_3$, we obtain the following

$$V_1 \geq (\frac{1}{8}\delta_b\delta_h - \frac{1}{2}k_1^{-2}\delta_h\Delta_h)\|X\|^2 + (\frac{1}{4}\delta_a\delta_b - \frac{1}{2}k_1^2\delta_a\Delta_a)\|Y\|^2 \geq 0 \qquad (3.5)$$

if we choose $\quad k_1^2 \leq \frac{1}{2}\frac{\delta_b}{\Delta_a}\quad$ and $\quad \Delta_h \leq \frac{1}{8}\frac{\delta_b^2}{\Delta_a}$,

in a similar way $\qquad V_2 \geq 0 \qquad (3.6)$

if we choose $\quad k_2^2 \leq \frac{1}{2}\delta_a \quad$ and $\quad \Delta_h \leq \frac{1}{16}\delta_a \delta_b$

so we have $\quad \Delta_h \leq k\delta_a \delta_b$

where $\quad k = \min\frac{1}{8}\left\{\frac{1}{2}, \frac{\delta_b}{\delta_a \Delta_a}\right\}, \quad (k<1), \quad$ and

$$V_3 \geq [\frac{1}{4}\delta_b \delta_h - (\Delta_b + 1)\sqrt{\varepsilon}]\|X\|^2 + [\frac{1}{4}\delta_a \delta_b - \frac{6\Delta_a + 7}{4}\sqrt{\varepsilon}]\|Y\|^2 + [\frac{1}{2}\delta_a - \sqrt{\varepsilon}]\|Z\|^2 \geq 0,$$

if we choose

$$\sqrt{\varepsilon} \leq \min\left\{\frac{\delta_b \delta_h}{4\Delta_b + 4}, \frac{\delta_a \delta_b}{6\Delta_a + 7}, \frac{\delta_a}{2}, 1\right\}.$$

Then, $\quad V_3 \geq \delta_4 (\|X\|^2 + \|Y\|^2 + \|Z\|^2)$ \hfill (3.7)

where, $\delta_4 = \min\left\{\frac{1}{4}\delta_b \delta_h - (\Delta_b + 1)\sqrt{\varepsilon}, \frac{1}{4}\delta_a \delta_b - \frac{6\Delta_a + 7}{4}\sqrt{\varepsilon}, \frac{1}{2}\delta_a - \sqrt{\varepsilon}\right\}.$

Finally, we are left with $V_4$. Since $P(t, X, Y, Z)$ satisfies $(2.2)$,
by Schwarz's inequality we obtain

$$\begin{aligned}|V_4| &\leq (\frac{1}{2}\Delta_b \|X\| + \Delta_a \|Y\| + 2\|Z\|)\|P(t, X, Y, Z)\| \\ &\leq \delta_5 (\|X\| + \|Y\| + \|Z\|)(\delta_0 + \delta_1(\|X\| + \|Y\| + \|Z\|)) \\ &\leq 3\delta_1 \delta_5 (\|X\|^2 + \|Y\|^2 + \|Z\|^2) + \sqrt{3}\delta_0 \delta_5 (\|X\|^2 + \|Y\|^2 + \|Z\|^2)^{\frac{1}{2}}\end{aligned}$$ \hfill (3.8)

where $\quad \delta_5 = \max\left\{\frac{1}{2}\Delta_b, \Delta_a, 2\right\}.$

Combining the inequalities $(3.5), (3.6), (3.7)$ and $(3.8)$ in $(3.4)$, we obtain

$$\dot{V} \leq -2\delta_6 (\|X\|^2 + \|Y\|^2 + \|Z\|^2) + \delta_7 (\|X\|^2 + \|Y\|^2 + \|Z\|^2)^{\frac{1}{2}}$$ \hfill (3.9)

where $\quad \delta_6 = \frac{1}{2}\min\{\delta_4, 3\delta_1 \delta_5\} \quad$ and $\quad \delta_7 = \sqrt{3}\delta_0 \delta_5.$

If we choose
$(\|X\|^2 + \|Y\|^2 + \|Z\|^2)^{\frac{1}{2}} \geq \delta_8 = 2\delta_7 \delta_6^{-1}$, inequality $(3.9)$ implies that

$$\dot{V} \leq -\delta_6 (\|X\|^2 + \|Y\|^2 + \|Z\|^2)$$ \hfill (3.10)

infact, we can obtain $\quad \dot{V} \leq -1 \quad$ if we choose

$$( \|X\|^2 + \|Y\|^2 + \|Z\|^2 )^{\frac{1}{2}} \geq \max\left\{\delta_6^{\frac{-1}{2}}, \delta_8\right\}.$$

Now we can prove that for any solution $V[X(t), Y(t), Z(t)]$ of (3.3) we ultimately have

$$( \|X\|^2 + \|Y\|^2 + \|Z\|^2 ) \leq \Delta_1$$

where $\Delta_1$ is a positive constant.

Suppose on the contrary, we would have $V(X(t), Y(t), Z(t)) \to \infty$, as $t \to \infty$, which contradicts inequality (3.2) that $V$ is non-negative. By using Yoshizawa's Theorem ([10] Theorem 15.8), we know that Eq. (1.1) has at least a periodic solution.

If $P(t, X, Y, Z) = P(t)$, let $[X_1(t), Y_1(t), Z_1(t)]$ and $[X_2(t), Y_2(t), Z_2(t)]$ be any solutions of (3.3), thus

$$\begin{cases} \dot{X}_1 = Y_1, \dot{Y}_1 = Z_1 \\ \dot{Z}_1 = -F(X_1, Y_1, Z_1)Z_1 - G(X_1, Y_1)Y_1 - H(X_1) + P(t) \end{cases}, \quad (3.11a)$$

$$\begin{cases} \dot{X}_2 = Y_2, \dot{Y}_2 = Z_2 \\ \dot{Z}_2 = -F(X_2, Y_2, Z_2)Z_2 - G(X_2, Y_2)Y_2 - H(X_2) + P(t) \end{cases} \quad (3.11b)$$

set $\psi = X_1 - X_2$, $\eta = Y_1 - Y_2$, $\tau = Z_1 - Z_2$, from (3.11) we obtain

$$\begin{cases} \dot{\psi} = \eta, \dot{\eta} = \tau \\ \dot{\tau} = -F(\psi, \eta, \tau)\tau - G(\psi, \eta)\eta - H(\psi) \end{cases} \quad (3.12)$$

**Remark:** Assume that Eq. (3.12) which obtained from Eq. (3.11) is true; because of the relevant literature. See Eq. (3.19) of [3], and, in particular $2V(\xi, \eta, \zeta)$ of [3; p. 268] and [9].

Then, rearranging the Lyapunov function in terms of $\psi, \eta, \tau$ we have

$$2V(\psi, \eta, \tau) = \frac{1}{4}\langle B\psi, B\psi \rangle + \frac{3}{2}\langle B\eta, \eta \rangle + \langle \tau, \tau \rangle + \langle \frac{1}{2}B\psi + A\eta + \tau, \frac{1}{2}B\psi + A\eta + \tau \rangle. \quad (3.13)$$

In view of (3.10) and (3.13) we have

$$\dot{V}(\psi, \eta, \tau) \leq -\delta V(\psi, \eta, \tau)$$

for some constant $\delta > 0$. By integrating both side of the inequality from $0$ to $t$

we obtain

$$V[\psi(t),\eta(t),\tau(t)] - V[\psi(0),\eta(0),\tau(0)] \leq -\delta\int_0^t V(\psi,\eta,\tau)dt$$

$$V[\psi(t),\eta(t),\tau(t)] \leq V[\psi(0),\eta(0),\tau(0)] - \delta\int_0^t V(\psi,\eta,\tau)dt$$

$$= K - \delta\int_0^t V(\psi,\eta,\tau)dt$$

and by using Gronwall-Reid Bellman inequality we can obtain

$$V[\psi(t),\eta(t),\tau(t)] \leq K\exp\left(-\delta\int_0^t V(\psi,\eta,\tau)dt\right)$$

$$\leq Ke^{-\delta t}.$$

Hence

$$\lim_{t\to\infty}\psi(t) = 0, \quad \lim_{t\to\infty}\eta(t) = 0, \quad \lim_{t\to\infty}\tau(t) = 0$$

and this is the required result.

From Lasalle's Theorem, we know that system (3.3) has a unique periodic solution.

The remaining of the proof can be completed by similar estimations arising in Tunç and Ateş [9].

## 4. Example.  For $n = 2$

$$F(X,Y,Z) = \begin{bmatrix} 2+x^2+y^2+z^2 & 0 \\ 0 & 2(2+x^2+y^2+z^2) \end{bmatrix}, \quad H(X) = \begin{bmatrix} x^2 \\ 2x^2 \end{bmatrix}$$

$$G(X,Y) = \begin{bmatrix} 1+x^2+y^2 & 0 \\ 0 & 2(1+x^2+y^2) \end{bmatrix}, \quad P(t,X,Y,Z) = \begin{bmatrix} xyz\cos(t+w) \\ 2xyz\cos(t+w) \end{bmatrix}$$

$\lambda_1(F) = 2+x^2+y^2+z^2 > 0$, $\lambda_2(F) = 2(2+x^2+y^2+z^2) > 0$,
$\lambda_1(G) = 1+x^2+y^2 > 0$, $\lambda_2(G) = 2(1+x^2+y^2) > 0$.

## Acknowledgment

The author thanks the referee for correcting errors in the original manuscript and for helpful suggestions.

# 5. References

[1] Feng, C. H., On the existence of almost periodic solutions of nonlinear third order differential equations. *Annals of Differential Equations*, 9 (1993), no.4, 420-424.

[2] Feng, C. H., The existence of periodic solutions for a third order nonlinear differential equations. *(Chinese )Gongcheng Shuxue Xuebao*,11 (1994), no.2, 113-117.

[3] Feng, C., On the existence of periodic solutions for a certain system of third order nonlinear differential equations. *Annals of Differential Equations*,11 (1995), no.3, 264-269.

[4] Mehri, B.; Shadman, D., Periodic solutions of certain third order nonlinear differential equations. *Studia Scientiarum Mathematicarum Hungarica*, 33 (1997), no.4, 345-350.

[5] Meng, F. W., Ultimate boundedness results for a certain system of third order nonlinear differential equations. *J. Math. Anal. Appl.* 177 (1993), no.2, 496-509.

[6] Tiryaki, A., Boundedness and periodicity results for a certain system of third order nonlinear differential equations. *Indian J. Pure Appl. Math.,* 30 (1999), no. 4, 361-372.

[7] Tunc. C., On the boundedness and periodicity of the solutions of a certain vector differential equation of third order. *Applied Mathematics and Mechanics* (English Edition) 20 (1999), no. 2, 163-170.

[8] Tunc, C.; Ates, M., Stability and boundedness results for solutions of certain third order nonlinear vector differential equations. *Nonlinear Dynamics, Springer Netherlands*, (2006), no. 45, 273-281.

[9] Tunc, C.; Ates, M., On the periodicity results for solutions of some certain third order nonlinear differential equations. *Advances in Mathematical Sciences and Applications*, 16 (2006), no.1, 1-14.

[10] Yoshizawa, T., Stability theory and the existence of periodic solutions and almost periodic solutions. *Applied Mathematical Sciences*, Vol.14. Springer-Verlag, NewYork-Heidelberg, 1975.

[11] Afuwape, A.U., Ultimate boundedness results for a certain systems of third- order onlinear differential equations. *Journal of Mathematical Analysis and Applications*, 97 (1983), no.1, 140-150.